\documentclass[12pt,reqno]{article}

\usepackage[usenames]{color}
\usepackage{amssymb}
\usepackage{amsmath}
\usepackage{amsthm}
\usepackage{amsfonts}
\usepackage{amscd}
\usepackage{graphicx}

\usepackage[colorlinks=true,
linkcolor=webgreen,
filecolor=webbrown,
citecolor=webgreen]{hyperref}

\definecolor{webgreen}{rgb}{0,.5,0}
\definecolor{webbrown}{rgb}{.6,0,0}

\usepackage{color}
\usepackage{fullpage}
\usepackage{float}

\usepackage{breakurl}

\newcommand{\Mod}[1]{\ \mathrm{mod}\ #1}
\begin{document}

\theoremstyle{plain}
\newtheorem{theorem}{Theorem}
\newtheorem{corollary}[theorem]{Corollary}
\newtheorem{lemma}[theorem]{Lemma}
\newtheorem{proposition}[theorem]{Proposition}
\theoremstyle{definition}
\newtheorem{definition}[theorem]{Definition}
\newtheorem{example}[theorem]{Example}
\newtheorem{conjecture}[theorem]{Conjecture}

\theoremstyle{remark}
\newtheorem{remark}[theorem]{Remark}

\begin{center}
\vskip 1cm{\LARGE\bf On a Curious Identity of Ramanujan
\vskip 1cm}
\large
H\`ung Vi\d{\^e}t Chu\\
Department of Mathematics\\
University of Illinois at Urbana-Champaign \\
Champaign, IL 61820\\
USA \\
\href{mailto:hungchu2@illinois.edu}{\tt hungchu2@illinois.edu} \\
\end{center}
\vskip .2 in

\begin{abstract}
Ramanujan wrote the following identity
\begin{align*}
\sqrt{2 \left(1 - \frac{1}{3^2}\right)
\left(1 - \frac{1}{7^2}\right)
\left(1 - \frac{1}{11^2}\right)
\left(1 - \frac{1}{19^2}\right)}
\  =
\
 \left(1 + \frac{1}{7}\right) \left(1 + \frac{1}{11}\right) \left(1 + \frac{1}{19}\right).
\end{align*} We find necessary and sufficient conditions for the integers in the identity and prove that there are only finitely many such identities, and provide a method to generate many interesting variations.
\end{abstract}

\section{Introduction}
In ranking mathematicians on the basis of pure talent, G. H. Hardy \cite{H} gave Ramanujan the highest score of 100. On this scale, he gave Hilbert a score of 80, himself a score of only 25, but said his colleague Littlewood merited 30. In fact, Ramanujan's work has fascinated generations of mathematicians even a century after his death.
The curious identity recorded in one of his notebooks \begin{align*}
\sqrt{2 \left(1 - \frac{1}{3^2}\right)
\left(1 - \frac{1}{7^2}\right)
\left(1 - \frac{1}{11^2}\right)
\left(1 - \frac{1}{19^2}\right)}
\  =
\
 \left(1 + \frac{1}{7}\right) \left(1 + \frac{1}{11}\right) \left(1 + \frac{1}{19}\right)
\end{align*} was mentioned by Berndt \cite{Be} who asked: ``Is this an isolated result, or are there other identities of this type?'' Reb\'ak \cite{Re} provided formulas that generate infinitely many similar identities and believed that the curious identity is related to the reciprocal of the Landau-Ramanujan constant
\begin{align*}
\frac{1}{K}\ =\ \sqrt{2}\prod_{\substack{p\mbox{ prime}\\p\equiv 3\Mod 4}}\sqrt{1-\frac{1}{p^2}}\ =\ & \sqrt{2\left(1-\frac{1}{3^2}\right)\left(1-\frac{1}{7^2}\right)\left(1-\frac{1}{11^2}\right)\left(1-\frac{1}{19^2}\right)\cdots}.
\end{align*}
Given $t, A, x, y, z\in\mathbb{R}$, we define
\begin{align*}
    f(t, A, x, y, z) &\ =\ \sqrt{t \left(1-\frac{1}{A^2}\right)\left(1-\frac{1}{x^2}\right)\left(1-\frac{1}{y^2}\right)\left(1-\frac{1}{z^2}\right)}, \mbox{ and }\\
    g(x,y,z)&\ =\ \left(1+\frac{1}{x}\right)\left(1+\frac{1}{y}\right)\left(1+\frac{1}{z}\right).
\end{align*}
In particular, Reb\'ak considered identities of the form
\begin{align}\label{mostgen}
f(t, A, x, y, z) \ =\ g(x,y,z).
\end{align}
However, Reb\'ak's generalization does not guarantee non-trivial integral values for all variables in Equation \eqref{mostgen}; in fact, there are only two integral solutions given by \cite[Theorem 2.1]{Re}. Motivated by this, we answer the question of whether or not there exist more integral solutions.

Our first main result states necessary and sufficient conditions for real numbers $t, A, x, y$, and $z$ that satisfy Equation \eqref{mostgen}. As we study necessary and sufficient conditions, our approach is more general and gives us a full understanding of the equation. It is worth noting that despite the complicated appearance of Equation \eqref{mostgen}, necessary and sufficient conditions for the variables can be established using quadratic equations. By a non-trivial solution to Equation \eqref{mostgen}, we mean that \begin{align}\label{nontri} A,x,y,z\notin\{0,\pm 1\}, t \neq 0 \mbox{ and }1+\frac{1}{z}  > 0.\end{align}
The first two conditions are reasonable. Let us explain the third condition. Observe that the first two conditions guarantee that a non-trivial solution makes both sides of Equation \eqref{mostgen} strictly positive. Because $g(x,y,z) > 0$ implies that at least one of its factors is positive,  we may assume that $1+1/z>0$ without loss of generality.

\begin{theorem}\label{sufnes}
The tuple $(t, A, x, y, z)$ is a non-trivial  solution to Equation \eqref{mostgen} if and only if
there exist $A,z\notin \{0,\pm 1\}$ with $1+1/z>0$ and $t,k\neq 0$ such that if we let
\begin{align*}
\gamma&\ =\ [(A^2-1)t-A^2]kz-[(A^2-1)t+A^2]k,\mbox{ and }\\
\beta&\ =\ [(A^2-1)t+A^2]kz-[(A^2-1)t-A^2]k-1,
\end{align*}
then
$$\begin{cases}\gamma^2-4\beta\ge 0,\\ (1+\gamma+\beta)\beta>0,
\end{cases}$$
and $x$ and $y$ are the roots of $X^2-\gamma X+\beta=0$.
\end{theorem}

\begin{remark}
We can use the sufficient condition of Theorem \ref{sufnes} to find many attractive identities. For example, with $t=1-\frac{1}{4^2}$, random searching\footnote{Programmed by G. Dresden at Washington and Lee University.} using ${\tt Mathematica}$ gives
\begin{align*}
\begin{split}
\sqrt{
\left(1 - \frac{1}{2^2}\right)
\left(1 - \frac{1}{3^2}\right)
\left(1 - \frac{1}{4^2}\right)
\left(1 - \frac{1}{9^2}\right)
\left(1 - \frac{1}{17^2}\right)}
\  = \
\left(1 - \frac{1}{3}\right)\left(1 + \frac{1}{9}\right) \left(1 + \frac{1}{17}\right).
\end{split}
\end{align*}\normalsize
\end{remark}

Next, define a \textit{perfect Ramanujan identity} to be an identity of Form \eqref{mostgen} with non-trivial positive integral values for all variables $A,x,y,z,t$. Define a \textit{general Ramanujan identity} to be an identity of Form \eqref{mostgen} with non-trivial integral values for all variables. We have the following theorem.

\begin{theorem}\label{finite}
There are only finitely many perfect Ramanujan identities, while there are infinitely many general Ramanujan identities.
\end{theorem}

The last section is devoted to some variations of the identity by Ramanujan, such as
\begin{align*}
\sqrt{\bigg(1-\frac{1}{9^2}\bigg)\bigg(1-\frac{1}{11^2}\bigg)\bigg(1-\frac{1}{23^2}\bigg)\bigg(1-\frac{1}{24^2}\bigg)\bigg(1-\frac{1}{45^2}\bigg)}\ = \ \bigg(1+\frac{1}{9}\bigg)\bigg(1-\frac{1}{11}\bigg)\bigg(1-\frac{1}{45}\bigg).\end{align*}
A theorem in the last section allows us to construct examples where the left side of the above identity is arbitrarily long. In the Appendix, we provide two curious identities and a comprehensive list of 39 super-perfect Ramanujan identities.

\section{Necessary and sufficient conditions}
\begin{proof}[Proof of Theorem \ref{sufnes}] Forward implication:
    from \eqref{mostgen}, we know that \begin{align*} t \bigg(1-\frac{1}{x}\bigg)\bigg(1-\frac{1}{y}\bigg)\bigg(1-\frac{1}{z}\bigg)\frac{A^2-1}{A^2}\ = \ \bigg(1+\frac{1}{x}\bigg)\bigg(1+\frac{1}{y}\bigg)\bigg(1+\frac{1}{z}\bigg).\end{align*}
So, \begin{align*} t\cdot \frac{A^2-1}{A^2}\cdot \frac{z-1}{z+1}\ = \ \frac{(1+\frac{1}{x})(1+\frac{1}{y})}{(1-\frac{1}{x})(1-\frac{1}{y})}.\end{align*}
Hence, \begin{align*} \frac{t(A^2-1)(z-1)}{A^2(z+1)}\ = \ \frac{xy+(x+y)+1}{xy-(x+y)+1}.\end{align*}
Therefore, there exists some non-zero $k\in\mathbb{R}$ such that
\begin{align*}
xy+(x+y)+1\ &=\ 2k t(A^2-1)(z-1),\\
xy-(x+y)+1\ &=\ 2kA^2(z+1).\end{align*} Simple computation shows that $xy=\beta$ and $x+y=\gamma$, so $x$ and $y$ are the roots of the equation $X^2-\gamma X+\beta=0$. Since $x$ and $y$ are real solutions to the quadratic equation $X^2-\gamma X+\beta= 0$, $\gamma^2-4\beta\ge 0$. Lastly, $1+1/z>0$ implies that $(1+1/x)(1+1/y)>0$. We may re-write this as
$$(1+1/x)(1+1/y) \ =\ \frac{(x+1)(y+1)}{xy} \ =\ \frac{\beta+\gamma+1}{\beta}\ >\ 0$$
to conclude that $(1+\gamma+\beta)\beta>0$. This completes our proof of the forward implication.

We now turn to the backward implication. For $z, A \notin\{0,\pm 1\}$ with $1+1/z>0$ and $t, k\neq 0$, define \begin{align*}\gamma&\ =\ [(A^2-1)t-A^2]kz-[(A^2-1)t+A^2]k, \mbox{ and }\\ \beta&\ =\ [(A^2-1)t+A^2)kz-[(A^2-1)t-A^2]k-1.\end{align*} By assumption,
    $$\begin{cases}\gamma^2-4\beta\ge 0,\\ (1+\gamma+\beta)\beta>0.
\end{cases}$$
    Let $x$ and $y$ be real roots of the quadratic equation $X^2-\gamma X+\beta=0$. Clearly, $x,y\notin \{0, \pm 1\}$. Then \begin{align*}
        x+y\ &= \ \gamma\ = \ [(A^2-1)t-A^2]kz-[(A^2-1)t+A^2]k,\\
        xy\ &= \ \beta\ = \ [(A^2-1)t+A^2)kz-[(A^2-1)t-A^2]k-1.
    \end{align*}
    So, \begin{align*}
        xy+(x+y)+1\ &= \ 2k t(A^2-1)(z-1),\\
        xy-(x+y)+1\ &= \ 2kA^2(z+1).
    \end{align*} We have
    \begin{align*}
    \frac{xy+(x+y)+1}{xy-(x+y)+1}\ = \ \frac{t(A^2-1)(z-1)}{A^2(z+1)}.
    \end{align*}
    Equivalently,
    \begin{align*}
    \frac{(x+1)(y+1)}{(x-1)(y-1)}\ = \ \frac{t(A^2-1)(z-1)}{A^2(z+1)}.
    \end{align*}
    Since $\beta\neq 0$, $x,y\neq 0$. We can safely divide the numerator and the denominator of the left side by $xy$ to have
    \begin{align*}
    \frac{(1+1/x)(1+1/y)}{(1-1/x)(1-1/y)}\ = \ \frac{t(A^2-1)(1-1/z)}{A^2(1+1/z)}.
    \end{align*}
    Therefore,
    \begin{align}\label{ue}
        \left(1+\frac{1}{x}\right)^2\left(1+\frac{1}{y}\right)^2\left(1+\frac{1}{z}\right)^2\ =\ t\left(1-\frac{1}{A^2}\right)\left(1-\frac{1}{x^2}\right)\left(1-\frac{1}{y^2}\right)\left(1-\frac{1}{z^2}\right).
    \end{align}
    Because
    $$\left(1+\frac{1}{x}\right)\left(1+\frac{1}{y}\right)\left(1+\frac{1}{z}\right) \ =\ \frac{1+\gamma+\beta}{\beta}\left(1+\frac{1}{z}\right) \ >\ 0,$$
    Equation \eqref{ue} implies Equation \eqref{mostgen}.
    This completes our proof.
\end{proof}

\begin{remark}
Reb\'ak's Theorem 2 \cite{Re} is a special case of our sufficient condition. If we let $A=a, t=\frac{a+1}{a-1}, z=6a+1$ and $k=\frac{1}{2a}$, then by Theorem \ref{sufnes}, $\gamma=5a+3$ and $\beta=6a^2+7a+2$. To satisfy the conditions in Theorem \ref{sufnes}, $a\in \left[\left(-\infty, \frac{2}{3}\right)\cup \left(-\frac{1}{2},-\frac{1}{3}\right)\cup\left(-\frac{1}{6},\infty\right)\right]\backslash \{0, \pm 1\}$. We obtain the identity
\begin{align*}
g\left(\frac{a+1}{a-1},a, 2a+1, 3a+2, 6a+1\right)\ =\ f(2a+1, 3a+2, 6a+1).
\end{align*}
\end{remark}

\begin{remark}
If we let $A=a,t=\frac{a+1}{a-1},z=6a+5$ and $k=\frac{1}{2a+2}$, then by Theorem \ref{sufnes}, $\gamma=5a+2$ and $\beta=6a^2+5a+1$. To satisfy the conditions in Theorem \ref{sufnes}, $$a\in \left[\left(-\infty, -1\right)\left(-\frac{2}{3},-\frac{1}{2}\right)\left(-\frac{1}{3},\infty\right)\right]\backslash\{0,1\}.$$ We obtain the identity
\begin{align*}
g\left(\frac{a+1}{a-1}, a, 2a+1, 3a+1, 6a+5\right) \ =\ f\left (2a+1, 3a+1, 6a+5\right).
\end{align*}
\end{remark}

\section{Finitely many perfect Ramanujan identities}

\subsection{Finitely many perfect Ramanujan identities}

\begin{proof}[Proof of Theorem \ref{finite}]
Without loss of generality, assume that $2\le x\le y\le z$.
We have
\begin{align*}
\begin{split}\label{perfectequation}
t&\ = \ \left(1+\frac{1}{A^2-1}\right)\left(1+\frac{2}{x-1}\right)\left(1+\frac{2}{y-1}\right)\left(1+\frac{2}{z-1}\right)\\
&\ \le \ \left(1+\frac{1}{A^2-1}\right)\left(1+\frac{2}{x-1}\right)^3.
\end{split}
\end{align*}
The first equality results from some elementary algebra to transform \eqref{mostgen} and the inequality results from $x\le y\le z$. Notice that $t\ge 2$.
We consider the two following cases.

Case 1: $A\le x$. Then
\begin{align*}
\begin{split}
2\ \le \ t&\ \le \ \left(1+\frac{1}{A^2-1}\right)\left(1+\frac{2}{x-1}\right)^3\\&\ \le\ \left(1+\frac{1}{A^2-1}\right)\left(1+\frac{2}{A-1}\right)^3.
\end{split}
\end{align*}
So, $A$ is bounded, which implies that there are finitely many values of $t$. Hence, there are finitely many values for $$\left(1+\frac{2}{x-1}\right)\left(1+\frac{2}{y-1}\right)\left(1+\frac{2}{z-1}\right).$$
Let $M$ be the finite set of possible values; for all $m\in M$, $m>1$. Fixing $m\in M$, we have
\begin{align*}
\begin{split}
1\ < \ m&\ = \ \left(1+\frac{2}{x-1}\right)\left(1+\frac{2}{y-1}\right)\left(1+\frac{2}{z-1}\right)\\&\ \le \ \left(1+\frac{2}{x-1}\right)^3.
\end{split}
\end{align*}
Thus there are finitely many values of $x$ corresponding to each $m$. Because $m\in M$, a finite set, there are finitely many values of $x$. Repeating the process, we can show that there are finitely many values of $y$ and $z$.

Case 2: $A>x$. We have
    \begin{align*}
    \begin{split}2\ \le  \ t&\ \le \ \left(1+\frac{1}{A^2-1}\right)\left(1+\frac{2}{x-1}\right)^3\\
    &\ \le \ \left(1+\frac{1}{x^2-1}\right)\left(1+\frac{2}{x-1}\right)^3 .
\end{split}
\end{align*}
Therefore $x$ is bounded, which implies there are finitely many values for $t$. Using the same argument as in Case 1, we can show that there are finitely many values of $A$, $y$ and $z$.

This completes our proof that there are finitely many perfect Ramanujan identities.

Next, we show that there are infinitely many general Ramanujan identities by simply giving a parameterization of a family of non-trivial solutions. Let $t=2,A=k,x=5,y=1-2k^2$ and $z=7$ for some $k\neq \pm 1$. The following is a correct identity
\begin{align*}
f(2, k, 5, 1-2k^2, 7)\ =\ g(5, 1-2k^2, 7).
\end{align*}
We have completed the proof of Theorem \ref{finite}.
\end{proof}

\subsection{Super-perfect Ramanujan identities}

A super-perfect Ramanujan identity requires stricter but reasonable conditions on the variables $1\le t<A<x<y<z$. The reason of the strict inequality is because if $x=y$, for example, $(1-\frac{1}{x^2})(1-\frac{1}{y^2})=(1-\frac{1}{x^2})^2$ and so, $(1-\frac{1}{x^2})$ can be put outside the square root. We find that there are exactly 39 super-perfect Ramanujan identities, which are provided in the Appendix.
\begin{remark}\label{tvalue}
We claim that for all super-perfect Ramanujan identities, $t=2$.
By simple manipulation, we know that if $(t,A,x,y,z)$ gives a super-perfect Ramanujan identity, then Equation \eqref{mostgen} is equivalent to
\begin{align*}\begin{split}
t&\ = \ \left(1+\frac{1}{A^2-1}\right)\left(1+\frac{2}{x-1}\right)\left(1+\frac{2}{y-1}\right)\left(1+\frac{2}
{z-1}\right)\\
&\ \le \ \left(1+\frac{1}{2^2-1}\right)\left(1+\frac{2}{3-1}\right)\left(1+\frac{2}{4-1}\right)\left(1+\frac{2}{5-1}\right)=\frac{20}{3}.\end{split}\end{align*}
On the other hand,
\begin{align*}\left(1+\frac{1}{A^2-1}\right)\left(1+\frac{2}{x-1}\right)\left(1+\frac{2}{y-1}\right)\left(1+\frac{2}
{z-1}\right)\ > \ 1.\end{align*}
Hence, $2\le t\le 6$. Consider $t\ge 3$. We know that $A\ge 4$. Then
\begin{align*}
\begin{split}
3\ \le \ t&\ = \ \left(1+\frac{1}{A^2-1}\right)\left(1+\frac{2}{x-1}\right)\left(1+\frac{2}{y-1}\right)\left(1+\frac{2}{z-1}\right)\\ &\ \le \left(1+\frac{1}{4^2-1}\right)\left(1+\frac{2}{5-1}\right)\left(1+\frac{2}{6-1}\right)\left(1+\frac{2}{7-1}\right)\ =\ 2.98\overline{6}\ <\ 3.
\end{split}
\end{align*}
This is a contradiction. Therefore, $t=2$.
\end{remark}

\section{Some variations}
In this section we provide several interesting variations of the Ramanujan identity. We first give an example of a variation and then present the main
theorem, which helps generate arbitrarily long identities.

\begin{remark}\label{goaldeterminant}
Let $\gamma$ and $\beta$ be defined as in Theorem \ref{sufnes}. We find that the discriminant of the equation $X^2-\gamma X+\beta=0$ is $$\gamma^2-4\beta\ = \ ((z-1)(A^2-1)kt-A^2k(1+z)-2)^2-8A^2k(1+z).$$
The goal is to cleverly choose values for $A,z,k$ and $t$ so that $\sqrt{\gamma^2-4\beta}$ is an integer in order to have integral values for $x$ and $y$.
\end{remark}

We illustrate Remark \ref{goaldeterminant} by an example. Let $A=a-1, z=2a+1, t=1-\frac{1}{a^2}$ and $k=\frac{1}{a^2-1}$. Substituting these values into our function for $\gamma^2-4\beta$, we have $\gamma^2-4\beta=16(2-a)$. Note that $\gamma=-2$ and $\beta=4a-7$. In order to satisfy the conditions in Theorem \ref{sufnes}, $$a\in \left[(-\infty, -1)\cup\left(-\frac{1}{2},\frac{7}{4}\right)\right]\backslash\{0,1\}.$$ Therefore, if we let $a=2-b^2$ for $b\in\mathbb{N}_{\ge 2}$, then $\sqrt{\gamma^2-4\beta}$ is a positive integer and $a$ satisfies all conditions. The following example gives an identity from this choice.

\begin{example}\label{firstexvari}
Let $b=5$ and so, $a=-23$. We have $z=2\cdot (-23)+1=-45, A=-23-1=-24$ and $t=1-\frac{1}{23^2}$. Also, $\gamma=-2$ and $\sqrt{\gamma^2-4\beta}=4\cdot 5=20$. By Theorem \ref{sufnes}, we know that $x=\frac{1}{2}(-2+20)=9$ and $y=\frac{1}{2}(-2-20)=-11$. Hence, we have the following identity.
\begin{align*}\begin{split}&\sqrt{\left(1-\frac{1}{9^2}\right)\left(1-\frac{1}{11^2}\right)\left(1-\frac{1}{23^2}\right)\left(1-\frac{1}{24^2}\right)\left(1-\frac{1}{45^2}\right)}\\& = \ \left(1+\frac{1}{9}\right)\left(1-\frac{1}{11}\right)\left(1-\frac{1}{45}\right).\end{split}\end{align*}
\end{example}
The following is a generalization of the above method. This theorem gives arbitrarily long identities.

\begin{theorem}\label{theoremvar}
Let $a=2-b^2$ for some positive integer $b\ge 2$. Let $n\in\mathbb{N}$ such that $a+1,a+2,...,a+n\neq 0$ and $a+n\neq 1$. Then the following identity holds
\begin{align*}
\begin{split}
&\sqrt{\bigg(1-\frac{1}{(2b+1)^2}\bigg)\bigg(1-\frac{1}{(2b-1)^2}\bigg)\bigg(1-\frac{1}{(2a+2n-1)^2}\bigg)\prod_{i=0}^n\bigg(1-\frac{1}{(a-1+i)^2}\bigg)}\\
&=\ \bigg(1-\frac{1}{2b+1}\bigg)\bigg(1+\frac{1}{2b-1}\bigg)\bigg(1+\frac{1}{2a+2n-1}\bigg).
\end{split}
\end{align*}
\end{theorem}

\begin{proof}
We use Theorem \ref{sufnes} to prove this identity. Let $$t\ =\ \bigg(1-\frac{1}{a^2}\bigg)\bigg(1-\frac{1}{(a+1)^2}\bigg)\cdots\bigg(1-\frac{1}{(a+n-1)^2}\bigg),$$ $A=a-1$, $z=2a+2n-1$ and $k=\frac{1}{(a-1)(a+n)}$. The conditions put on $a$ and $n$ by Theorem \ref{theoremvar} satisfy the conditions mentioned in Theorem \ref{sufnes}. Plugging these values for $t,a,z$, and $k$ into our formula for $\gamma^2-4\beta$, we have $\sqrt{\gamma^2-4\beta}=\sqrt{16(2-a)}=\sqrt{16b^2}=4b$ and $\gamma=-2$. Therefore, $x=\frac{1}{2}(-2-4b)=-2b-1$ and $y=\frac{1}{2}(-2+4b)=2b-1$. So, by Theorem \ref{sufnes}, we have the identity.
\end{proof}

\section{Conclusion}

We end by discussing two questions suggested by the anonymous referee.
\begin{enumerate}
    \item For $n\ge 3$, define
    \begin{align*}f(t, A, x_1, x_2, \ldots, x_n) &\ =\ \sqrt{t\left(1-\frac{1}{A^2}\right)\left(1-\frac{1}{x_1^2}\right)\cdots \left(1-\frac{1}{x_n^2}\right)},\\
    g(x_1, x_2, \ldots, x_n) &\ =\ \left(1+\frac{1}{x_1}\right)\left(1+\frac{1}{x_2}\right)\cdots \left(1+\frac{1}{x_n}\right).
    \end{align*}
    Let $S$ be the set of all $(t, A, x_1, x_2, \ldots, x_n)$ that satisfies
    $$\begin{cases}
    f(t, A, x_1, x_2, \ldots, x_n) = g(x_1, x_2, \ldots, x_n),\\
    t\ge 1,\\
    A, x_1,  x_2, \cdots, x_n\ge 2.
    \end{cases}$$
    Is $S$ still finite when $n\ge 4$? Using our proof of Theorem \ref{finite} and induction, we know that $S$ is finite for all $n\ge 4$.
    \item What are identities of the form
    $$\sqrt{t\prod_{k=1}^n\left(1-\frac{1}{x_k^2}\right)}\ =\ \prod_{k=1}^n\left(1+\frac{1}{x_k}\right)?$$
    The author attempted this problem with $n = 4$ by introducing polynomials of higher degree but was unable to characterize all solutions as in the case $n = 3$. Finding some special families of solutions is possible with telescoping products and has been done by Reb\'ak \cite[Theorem 3 and Theorem 4]{Re}.
\end{enumerate}

\section{Appendix}

\subsection{Two identities inspired by Ramanujan}
Below are two curious identities. Since proving them is not the focus of the paper and the proof is itself not interesting, we simply mention them here.

\begin{theorem}
For $a\ge 3$, we have
\begin{align*}
    \begin{split}
&\sqrt{\bigg(1-\frac{1}{a^2}\bigg)\bigg(1-\frac{1}{(a-1)^2}\bigg)\bigg(1-\frac{1}{(2a+1)^2}\bigg)\bigg(1-\frac{1}{(1+2\sqrt{a-1})^2}\bigg)\bigg(1-\frac{1}{(2\sqrt{a-1}-1)^2}\bigg)}\\
&=\ \bigg(1+\frac{1}{2a+1}\bigg)\bigg(1+\frac{1}{1+2\sqrt{a-1}}\bigg)
\bigg(1-\frac{1}{2\sqrt{a-1}-1}\bigg).
\end{split}
\end{align*}\normalsize
\end{theorem}
\begin{theorem}
For $a\le 1$, we have
\begin{align*}
\begin{split}
    &\sqrt{\bigg(1-\frac{1}{a^2}\bigg)\bigg(1-\frac{1}{(a-1)^2}\bigg)\bigg(1-\frac{1}{(2a+1)^2}\bigg)\bigg(1-\frac{1}{(2\sqrt{2-a}+1)^2}\bigg)\bigg(1-\frac{1}{(2\sqrt{2-a}-1)^2}\bigg)}\\
    &=\
    \bigg(1-\frac{1}{2\sqrt{2-a}+1}\bigg)\bigg(1+\frac{1}{2\sqrt{2-a}-1}\bigg)
    \bigg(1+\frac{1}{2a+1}\bigg).
\end{split}
\end{align*}\normalsize
\end{theorem}

\subsection{List of super-Perfect Ramanujan identities}
We provide a comprehensive list of super-perfect Ramanujan-type identities. There are 39 of them. To find these, we use Remark \ref{tvalue} and the proof of Theorem \ref{finite} to find upper and lower bounds for all variables. If we require that all variables are prime like the original one from Ramanujan, there are only three of them (the last three).
\begin{center}
\begin{tabular}{ |c|c|c|c|c||c|c|c|c|c| }
\hline
$t$ & $A$ & $x$ & $y$ & $z$ & $t$ & $A$ & $x$ & $y$ & $z$\\
 \hline
 2 & 6 & 7 & 9 & 13 & 2 & 3 & 4 & 39 & 151\\
 2 & 5 & 6 & 7 & 71 & 2 & 3 & 4 & 41 & 127\\
 2 & 5 & 6 & 8 & 31 & 2 & 3 & 4 & 43 & 111\\
 2 & 5 & 6 & 11 & 15 & 2 & 3 & 4 & 46 & 95\\
 2 & 4 & 5 & 10 & 89 &  2& 3 & 4 & 47 & 91\\
 2 & 4 & 5 & 11 & 49 & 2& 3 & 4 & 51 & 79\\
 2 & 4 & 5 & 13 & 29 & 2& 3 & 4 & 55 & 71\\
 2 & 4 & 5 & 14 & 25 & 2& 3 & 4 & 61 & 63\\
 2 & 4 & 5 & 17 & 19 & 2& 3 &5 & 12 & 703\\
 2 & 4 & 6 & 7 & 449 & 2& 3 & 5 & 15 & 55\\
 2 & 4 & 6 & 8 & 49 & 2& 3 & 6 & 9 & 127\\
 2 & 4 & 6 & 9 & 29 & 2& 3 & 6 & 15 & 19\\
 2 & 4 & 7 & 9 & 17 & 2& 3 & 7 & 8 & 55\\
 2 & 4 & 8 & 9 & 13 & 2& 3 & 7 & 9 & 31\\
 2 & 3 & 4 & 32 & 991 & 2 & 3 & 7 & 10 & 23\\
 2 & 3 & 4 & 33 & 511 & 2 & 3 & 7 & 13 & 15\\
 2 & 3 & 4 & 34 & 351 & 2 & 3 & 5 & 13 & 127\\
 2 & 3 & 4 & 35 & 271 & 2 & 3 & 5 & 19 & 31\\
 2 & 3 & 4 & 36 & 223 &  2 & 3 & 7 & 11 & 19\\
 2 & 3 & 4 & 37 & 191 & & & & &\\
 \hline
\end{tabular}
\end{center}

\section{Acknowledgments}
The author wishes to thank Prof. Aaron Abrams, Prof. Kevin Beanland and Prof. Gregory Dresden at Washington and Lee University for many helpful conversations that raise new ideas for this paper. Thanks to Prof. Steven Miller at Williams College for comments on an earlier draft. Finally, thanks to Prof. Reb\'ak for pointing an error in the Arxiv version and the anonymous referee for useful suggestions.

\bigskip
\hrule
\bigskip

\noindent 2010 {\it Mathematics Subject Classification}: Primary 11A67; Secondary 11E25.

\noindent \emph{Keywords:} Ramanujan identity, reciprocal, square root

\bigskip
\hrule
\bigskip

\end{document}